\documentclass{article}
\usepackage{amsmath,amsthm,amsfonts}
\usepackage{enumerate}
\usepackage{amssymb}
\usepackage{epsf,psfig}

\newcommand{\R}{{\mathbb{R}}}
\newcommand{\N}{{\mathbb{N}}}

\newtheorem{assumption}{Assumption}
\newtheorem{proposition}{Proposition}
\newtheorem{definition}{Definition}

\newtheorem{theorem}{Theorem}
\newtheorem{lemma}{Lemma}
\newtheorem{claim}{Claim}
\newtheorem{remark}{Remark}
\newtheorem{data}{Data}
\newtheorem{example}{Example}

\newcommand{\LL}{{\mathcal{L}}}
\newcommand{\FF}{{\mathcal{F}}}

\newcommand{\nodes}{{\mathcal{N}}}
\newcommand{\arcs}{{\mathcal{A}}}
\newcommand{\neighbors}{{\mathrm{Neighbors}}}
\newcommand{\convexhull}{{\mathrm{conv}}}
\newcommand{\setarrow}{\rightrightarrows}
\newcommand{\closedsubsets}{{\mathcal{P}}}
\newcommand{\image}{{\mathrm{Image}}}
\newcommand{\ball}{{\mathcal{B}}}
\newcommand{\mass}{{\mathcal{X}}}
\newcommand{\sass}{X}
\newcommand{\diam}{\mu}
\newcommand{\ri}{{\mathrm{ri}}}

\newcommand{\GG}{{\mathcal{G}}}

\begin{document}

\title{Time-dependent unidirectional communication in multi-agent
systems} 

\author{Luc Moreau\footnote{Postdoctoral Fellow of the Fund for
Scientific Research - Flanders (Belgium) (F.W.O.-Vlaanderen). This
paper presents research results of the Belgian Programme on
Inter-University Poles of Attraction, initiated by the Belgian State,
Prime Minister's Office for Science, Technology and Culture. The
scientific responsibility rests with its authors.  EESA-SYSTeMS, Ghent
University, Technologiepark~914, 9052~Zwijnaarde, Belgium
(Tel: +32.9.264.56.55; Fax: +32.9.264.58.40; Email: Luc.Moreau@UGent.be).}\\[2ex]
Ghent University, Belgium}
\maketitle

\begin{abstract}
We study a simple but compelling model of $n$~interacting agents via
time-dependent, unidirectional communication.  The model finds wide
application in a variety of fields including synchronization, swarming
and distributed decision making.  In the model, each agent updates his
current state based upon the current information received from other
agents.  Necessary and/or sufficient conditions for the convergence of
the individual agents' states to a common value are presented,
extending recent results reported in the literature.  Unlike previous,
related studies, the approach of the present paper does not rely on
algebraic graph theory and is of a completely nonlinear nature. It is
rather surprising that with these nonlinear tools, extensions may be
obtained even for the linear cases discussed in the literature.  The
proof technique consists of a blend of graph-theoretic and
system-theoretic tools integrated within a formal framework of
set-valued Lyapunov theory and may be of independent interest.
Further, it is also observed that more communication does not
necessarily lead to better convergence and may eventually even lead to
a loss of convergence, even for the simple models discussed in the
present paper.  \\[2ex] {\bf{Keywords:\/}} multi-agent systems,
stability analysis, swarms, synchronization, set-valued Lyapunov
theory.
\end{abstract}

\section{Introduction}

Recent years have witnessed an increasing interest in the interaction
between information flow and system dynamics.  It is recognized that
information and communication constraints may have a considerable
impact on the performance of a control system. The study of these
topics forms a very active area of research, giving rise to new
control paradigms such as {\em quantized control systems}, {\em
networked control systems} and {\em multi-agent (multi-vehicle,
multi-robot) systems}.

An important aspect of information flow in a dynamical system is the
communication topology, which determines what information is available
for which component at a given time instant. We study the role of
communication topology in the context of multi-agent systems. Our
interest in multi-agent systems is motivated by the emergence of
several applications including formation flying of UAVs (unmanned
aerial vehicles),
cooperative robotics~\cite{CaFuKa:97,SuKu:99,FoBuKrTh:00} and sensor
networks~\cite{asp:03}.  

In the present paper we consider a group of $n$~agents, not
necessarily identical. The individual agents share a common state
space and each agent updates his current state based upon the
information received from other agents, according to a simple rule.
Our aim is to relate the information flow and communication
structure with the stability properties of the group of agents.

The model that we consider encompasses, or is closely related to,
several models reported in the literature.
A prominent and well-studied example concerns synchronization of
coupled oscillators, a phenomenon which is ubiquitous in the natural
world and finds several applications in physics and engineering;
see~\cite{St:00,St:03} for a review and a historical account. The
Kuramoto equation~\cite{Ku:84,St:00}, which is widely accepted in that
field, models a population of oscillators with (typically sinusoidal)
coupling terms.  The Kuramoto equation may be studied with the tools
of the present paper.  Another interesting example concerns swarming,
a cooperative behavior observed for a variety of living beings such as
birds, fish, bacteria, etc.  In the physics literature, swarming
models are often individual-based with each individual being
represented by a particle moving with constant velocity, its direction
of motion being updated according to nearest neighbor coupling; see,
for example, the paper~\cite{ViCzJaCoSc:95}.  In recent years,
engineering applications such as formation control have increased the
interest of engineers in swarms~\cite{LeFi:01,JaLiMo:03,GaPa:03}.  The
approach of the present paper applies to the swarming models discussed
in~\cite{ViCzJaCoSc:95,JaLiMo:03}.  
A third and last example of
multi-agent systems that fall within the scope of the present paper
are the concensus algorithms studied in~\cite{OlMu:acc03}.  Consensus
protocols enable a network of dynamic agents to agree upon quantities
of interest via a process of distributed decision making.  Agreement
problems have a long history in the field of computer science and are
currently finding applications in formation control of multi-vehicle
systems~\cite{FaMu:02a,FaMu:02b}.  

We impose very weak assumptions on the communication topology; we
allow for unidirectional and time-dependent communication patterns.
Unidirectional communication is important for practical applications
and can easily be incorporated, for example, via broadcasting.  Also,
sensed information flow which plays a central role in schooling and
flocking, is typically not bidirectional.  In addition, we do not
exclude loops in the communication topology.  This means that,
typically, we are considering leaderless coordination rather than a
leader-follower approach.  (Leaderless coordination is also
considered, for example, in~\cite{FaMu:02a,OlMu:ifac02,BaLe:02}).
Finally, we allow for time-dependent communication patterns which are
important if we want to take into account link failure and link
creation, reconfigurable networks and nearest neighbor coupling.

The contribution of the paper is twofold.  We present necessary and/or
sufficient conditions on the communication topology guaranteeing
convergence of the individual agents' states to a common value. The
results that we obtain extend the results reported
in~\cite{Ku:84,St:00,JaLiMo:03,OlMu:acc03}.  A second contribution
concerns the proof technique enabling us to obtain these results.
Previous studies that have
considered the connection between communication topology and system
stability such as~\cite{FaMu:02a,JaLiMo:03,OlMu:acc03} typically rely on algebraic graph theory, relating the graph
topology with the algebraic structure of associated graph matrices.
In the present paper, instead of relying on algebraic graph theory, we
propose a blend of graph-theoretic and system-theoretic tools to
analyse stability.  Our approach is of an inherently nonlinear nature
with the notion of convexity playing a central role.  It is rather
surprising that with these nonlinear tools, extensions may be obtained
even for the linear cases discussed in the literature.  The analysis
tools that we propose are integrated within a formal framework relying
on set-valued Lyapunov functions and may be of independent interest.

The paper is organized as follows.  Section~\ref{s:linear} illustrates
the main results of the paper by means of a simple example.
Section~\ref{s:model} introduces the model that is studied in the
remainder of the paper.  Section~\ref{s:stab} proposes a novel
stability concept and a corresponding Lyapunov characterization; this
may be of independent interest.
Sections~\ref{s:uga}--\ref{s:bidirectional} present necessary and/or
sufficient conditions for convergence of the agents' states to a
common, constant value.  Section~\ref{s:conclusion} summarizes the
main results of the paper.  Appendix~\ref{s:appendix} contains a
graph-theoretic result that is used in the technical proofs of the
paper.

\subsection*{Notation and terminology}

We distinguish between the set-inclusion symbols~$\subset$
and~$\subseteq$. The symbol~$\subset$ denotes strict inclusion
whereas $\subseteq$ denotes non-strict inclusion.

Let~$K$ be a subset of a finite-dimensional Euclidean space~$X$.  The
set~$\ball(K,c)$ ($c>0$) is defined as the set of points in~$X$ whose
distance to~$K$ is strictly smaller than~$c$. When~$K$ is a
singleton~$\{x\}$ the notation~$\ball(x,c)$ is used instead
of~$\ball(\{x\},c)$.

Consider finite-dimensional Euclidean spaces~$X$ and~$Y$.  We denote
by~$\closedsubsets(Y)$ the collection of all closed subsets of~$Y$.  A
{\em{set-valued function}}~$f:X\setarrow Y$ is a (single-valued)
function~$f:X\rightarrow\closedsubsets(Y)$.  A set-valued
function~$f:X\setarrow Y$ is called {\em{upper semicontinuous}} if for
every~$x\in X$ and every~$\varepsilon>0$ there is~$\delta>0$ such that
$f(y)\in\ball(f(x),\varepsilon)$ whenever $y\in\ball(x,\delta)$.

The {\em relative interior} of a convex set~$C$ in finite-dimensional
Euclidean space is the interior which results when this set is
regarded as a subset of its affine hull and is denoted by~$\ri(C)$. A
{\em polytope} is a set which is the convex hull of finitely many
points. See~\cite{Ro:97} for more information.

\section{Example}
\label{s:linear}

The present section illustrates the main results of the paper by means
of a linear example.  The results that we obtain for this example
generalize results that have recently been reported in the
literature.
At the same time, the present section
also introduces several notions that will be used repeatedly in the
remainder of the paper (directed graph, neighbors, weak connectivity).
\begin{definition}[directed graph]
A directed graph is a pair~$(\nodes,\arcs)$ where
$\nodes$ is a nonempty, finite set and
$\arcs$ is a subset of $\nodes\times\nodes$ satisfying
$(k,k)\notin\arcs$ for all~$k\in\nodes$.
Elements of~$\nodes$ are referred to as nodes and an element~$(k,l)$ of~$\arcs$
is referred to as an arc from~$k$ to~$l$.
\end{definition}\noindent
\begin{definition}[neighbors]
Consider a directed graph~$(\nodes,\arcs)$ and a nonempty
subset~$\LL\subseteq\nodes$. The set~$\neighbors(\LL,\arcs)$ is the set
of those nodes~$k\in\nodes\setminus\LL$ for which there is~$l\in\LL$
such that $(k,l)\in\arcs$.  
When~$\LL$ is a singleton~$\{m\}$, the
notation~$\neighbors(m,\arcs)$ is used instead
of~$\neighbors(\{m\},\arcs)$.
\end{definition}\noindent
Notice that any
nonempty subset~$\LL\subseteq\nodes$ satisfies
\[
\neighbors(\LL,\arcs)
=
\left(\cup_{m\in\LL}\neighbors(m,\arcs)\right)\setminus\LL.
\]
\begin{definition}[weighted directed graph]
A weighted directed graph is a triple~$(\nodes,\arcs,w)$ where
$(\nodes,\arcs)$ is a directed graph and~$w:\arcs\rightarrow\R_{>0}$
is a map associating to each arc~$(k,l)$ a positive
weight denoted by~$w_{kl}$.
\end{definition}\noindent

\begin{figure}[h]
\hspace*{\fill}
\psfig{file=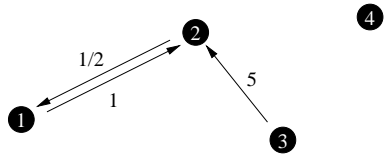,height=2cm}
\hspace*{\fill}
\caption{Weighted directed graph}
\label{f:weighteddirected}
\end{figure}

Consider, for example, the weighted directed graph
depicted in Fig.~\ref{f:weighteddirected}.
Each node in this graph corresponds to an agent and each arc
represents a communication channel.  Each agent~$i$ ($i=1,\dots,4$) is being
attributed a real state variable~$x_i\in\R$.  Suppose that this
weighted directed graph represents the communication pattern at
time~$t$.  Agent~$1$ receives information from agent~$2$ and updates
his state~$x_1$ according to the weighted average
\[
x_1(t+1)=\frac{1x_1(t)+(1/2)x_2(t)}{1+1/2}.
\]
Notice that agent~$1$'s own state~$x_1(t)$ is taken into account with unit
weight~$1$.  Agent~$2$ receives information from agents~$1$ and~$3$
and updates his heading~$x_2$ according to the weighted average
\[
x_2(t+1)=\frac{1x_1(t)+1x_2(t)+5x_3(t)}{1+1+5}.
\]
Agents~$3$ and~$4$ do not receive any information:
\begin{align*}
x_3(t+1)&=x_3(t),\\
x_4(t+1)&=x_4(t).
\end{align*}
We are interested in the dynamics of the $x_i$-variables of the agents
that arises when the communication graph changes over time.  In order
to formulate the dynamics for general graphs, we
associate to a weighted directed graph~$\GG=(\nodes,\arcs,w)$ with
vertex set~$\nodes=\{1,\dots,n\}$ a real $n\times n$-matrix~$A(\GG)$
whose components we define as follows:
\begin{equation}
(A(\GG))_{kl}=
\begin{cases}
\frac{w_{lk}}{1+\sum_{i\in\neighbors(k,\arcs)}w_{ik}}
&\mbox{whenever $(l,k)\in\arcs$,}\\
\frac{1}{1+\sum_{i\in\neighbors(k,\arcs)}w_{ik}}
&\mbox{whenever $k=l$,}\\
0&\mbox{otherwise.}
\end{cases}
\label{e:a}
\end{equation}
Notice that the matrix~$A(\GG)$ is {\it{stochastic\/}};
that is, it is a square matrix with non-negative entries and with the
property that its row sums are all equal to~$1$.  For the weighted
directed graph of Fig.~\ref{f:weighteddirected}, for example, the
matrix~$A(\GG)$ becomes
\begin{equation*}
\begin{pmatrix}
2/3&1/3&0&0\\
1/7&1/7&5/7&0\\
0&0&1&0\\
0&0&0&1
\end{pmatrix}.
\end{equation*}

Consider a sequence of weighted directed
graphs~$\GG(t)=(\nodes,\arcs(t),w(t))$ with common vertex
set~$\nodes=\{1,\dots,n\}$ and where~$t\in\N$.  We are interested in
the following discrete-time system on~$\R^n$:
\begin{equation}
x(t+1)=A(\GG(t))x(t).\\
\label{e:linearconsensus}
\end{equation}
In particular, we want to formulate necessary and/or sufficient
conditions for the convergence of the $n$~components~$x_i(t)$ to a
common value as~$t\rightarrow\infty$.  It is clear that the
convergence properties of Eq.~\eqref{e:linearconsensus} depend both on
the connectivity properties of the sequence of directed graphs and on
the associated weight functions. In the present paper we focus on the
role of the connectivity properties and assume that the
weights are well-behaved by imposing a uniform boundedness condition.
Below we state two results that follow as a consequence of the general
theory established in the paper (Theorems~\ref{t:NS-UGA}
and~\ref{t:bidirectional}). The formulation of these results relies on
the notion of weak connectivity.
\begin{definition}[weak connectivity]
Consider a directed graph~$(\nodes,\arcs)$.  A node~$k\in\nodes$ is
connected to a node~$l\in\nodes$ if there is a path from~$k$ to~$l$ in
the graph which respects the orientation of the arcs.
A directed graph~$(\nodes,\arcs)$ is called weakly connected if there
is a node~$k\in\nodes$ which is connected to all other
nodes~$l\in\nodes\setminus\{k\}$. A sequence of directed
graphs~$(\nodes,\arcs(t))$ with~$t\in\N$ is called weakly connected
across an interval~$I\subseteq\N$ if the directed
graph~$(\nodes,\cup_{t\in I}\arcs(t))$ is weakly connected.
\end{definition}\noindent
\begin{proposition}[unidirectional]\label{p:lc1}
Consider a sequence of weighted directed
graphs~$\GG(t)=(\nodes,\arcs(t),w(t))$ with common vertex
set~$\nodes=\{1,\dots,n\}$ and where~$t\in\N$.  Assume the existence
of~$0<e_\mathrm{min}\leq e_\mathrm{max}<\infty$ such that the weight
functions~$w(t)$ take values in~$[e_\mathrm{min},\,e_\mathrm{max}]$
for all~$t\in\N$.
If there is $T\geq 0$ such that for all $t_0\in\N$ the sequence of
directed graphs~$(\nodes,\arcs(t))$ is weakly connected
across~$[t_0,\,t_0+T]$, then the $n$-components~$\zeta_i(t)$ of any
solution~$\zeta$ of Eq.~\eqref{e:linearconsensus} converge to a common
value as~$t\rightarrow\infty$.
\end{proposition}\noindent
\begin{proposition}[bidirectional]\label{p:lc2}
Consider a sequence of weighted directed
graphs~$\GG(t)=(\nodes,\arcs(t),w(t))$ with common vertex
set~$\nodes=\{1,\dots,n\}$ and where~$t\in\N$.  Assume the existence
of~$0<e_\mathrm{min}\leq e_\mathrm{max}<\infty$ such that the weight
functions~$w(t)$ take values in~$[e_\mathrm{min},\,e_\mathrm{max}]$
for all~$t\in\N$.
Assume in addition that the graphs~$(\nodes,\arcs(t))$ are
bidirectional for all~$t\in\N$.  If for all $t_0\in\N$ the sequence of
graphs~$(\nodes,\arcs(t))$ is weakly connected
across~$[t_0,\,\infty)$, then the $n$-components~$\zeta_i(t)$ of any
solution~$\zeta$ of Eq.~\eqref{e:linearconsensus} converge to a common
value as~$t\rightarrow\infty$.
\end{proposition}\noindent
\subsection{Discussion of Propositions~\ref{p:lc1} and~\ref{p:lc2}}
\begin{enumerate}
\item
Let us mention some of the subtleties involved in the study of
Eq.~\eqref{e:linearconsensus}. First of all,
Eq.~\eqref{e:linearconsensus} belongs to the class of linear
time-varying systems, whose stability properties are hard to analyse
in general.  Secondly, the equilibrium value to which the components
are shown to converge, is not known explicitly; this equilibrium value
depends both on the initial data and on the sequence of
weighted directed graphs.  Thirdly, it has been shown in~\cite{JaLiMo:03}
that, in general, there does not exist a time-invariant, quadratic
Lyapunov function for Eq.~\eqref{e:linearconsensus}.  Partially
motivated by this negative result, we put forward in the present paper
a series of graph-theoretic and system-theoretic analysis tools which
are {\em{not}} based on the linear structure of
Eq.~\eqref{e:linearconsensus}. This contrasts our approach with the
approaches taken in~\cite{JaLiMo:03,OlMu:acc03} which rely
heavily on linear algebra and algebraic graph theory.
\item
For the special case that the weights are all equal to one,
Eq.~\eqref{e:linearconsensus} reduces to the model studied
in~\cite[Theorems~1--2]{JaLiMo:03}.  
From a control perspective, an important contribution
of~\cite{JaLiMo:03} is that this work relates the information flow
and communication structure with the stability properties of the group
of agents. 
The results obtained in the
present paper generalize the results
from~\cite[Theorems~1--2]{JaLiMo:03} in several directions.  Apart
from the observation that the present models are more general, it is
important to notice that unlike~\cite{JaLiMo:03} we do not assume
bidirectional communication (Proposition~\ref{p:lc1}) and, when we
consider bidirectional communication, we do not assume uniformity in
the connectivity condition (Proposition~\ref{p:lc2}).  It is rather
surprising that with our inherently nonlinear approach, we are able to
obtain generalizations of results recently reported in the literature,
even when restricting attention to the linear case.
\item
As a corollary of Proposition~\ref{p:lc1} we obtain the following
graph theoretic result.  A stochastic $n\times n$-matrix~$A$ with
positive diagonal elements has all but one of its eigenvalues strictly
inside the unit circle (the only exception being the trivial
eigenvalue~$1$) if (and only if---see~Theorem~\ref{t:NS-UGA}) the
associated directed graph~$(\nodes,\arcs)$
\[
\nodes=\{1,\dots,n\},\qquad\arcs=\{(k,l)\in\nodes\times\nodes:A_{lk}>0\mbox{
and }k\not =l\}
\]
is weakly connected.  As before, we mention that we obtain this graph
theoretic result using inherently nonlinear analysis tools. Moreover,
the present approach enables us to go beyond this result by
considering time-dependent stochastic matrices.
\end{enumerate}

\section{Multi-agent dynamics}
\label{s:model}

We introduce here the dynamical equations that will be studied in the
remainder of the paper.  These equations model the interaction of
$n$~agents via time-dependent, unidirectional communication.  
The individual agents share a common state space~$\sass$ which is
assumed to be finite-dimensional Euclidean. We are interested in
multi-agent dynamics that guarantee convergence of the individual
agents' states to a common value.  The dynamics considered in the
present paper are governed by a discrete-time system on~$\sass^n$. Its
formulation involves two ingredients: a time-sequence of directed
graphs and a discrete-time map~$f$ determining the actual dynamics.
\begin{data}[communication graphs]\label{d:a}
A sequence of directed graphs~$(\nodes,\arcs(t))$ with common node
set~$\nodes=\{1,\dots,n\}$ and with $t\in\N$;
\end{data}\noindent
\begin{data}[discrete-time map]\label{d:f}
A continuous function~$f:\N\times\sass^n\rightarrow\sass^n$.
\end{data}\noindent
These ingredients give rise to the following discrete-time system
on~$\sass^n$
\begin{equation}
x(t+1)=f(t,x(t))\label{e:consensus}
\end{equation}
or, expressed in terms of the individual agents' states,
\begin{equation*}
\begin{split}
x_1(t+1)&=f_1(t,x_1(t),\dots,x_n(t)),\\
        &\vdots\\
x_n(t+1)&=f_n(t,x_1(t),\dots,x_n(t)).
\end{split}
\end{equation*}
Here we have introduced the decomposition $x=(x_1,\dots,x_n)$ and
$f=(f_1,\dots,f_n)$ corresponding to the product structure
of~$\sass^n$.

The directed
graph~$(\nodes,\arcs(t))$ characterizes the inter-agent communication
at time~$t$. Each node of this graph corresponds to an agent and each
arc represents a communication channel.  The communication
graph~$(\nodes,\arcs(t))$ determines the information that is available
for the agents at time~$t$ and each agent updates its state based upon
this information according to Eq.~\eqref{e:consensus}. 
We impose two assumptions on the
map~$f=(f_1,\dots,f_n)$.
\begin{assumption}[communication]\label{a:f1}
For every time~$t\in\N$ and for each agent~$k\in\{1,\dots,n\}$
the value of the function~$f_k(t,\cdot):\sass^n\rightarrow \sass$ depends only
on the states of agent~$k$ and agents~$i\in\neighbors(k,\arcs(t))$; 
%
that is, 
$f_k(t,x)=f_k(t,\overline{x})$ whenever $x_k=\overline{x}_k$ and
$x_i=\overline{x}_i$ for all $i\in\neighbors(k,\arcs(t))$.
%
\end{assumption}\noindent
\begin{assumption}[convexity]\label{a:f2}
Associated to each directed
graph~$(\nodes,\arcs)$ with node set~$\nodes=\{1,\dots,n\}$ and each
agent~$k$ there is a 
continuous set-valued
function~$e_k(\arcs):\sass^n\setarrow\sass$ satisfying
\begin{enumerate}
\item
\begin{equation}
f_k(t,x)
\in e_k(\arcs(t))(x),\qquad\forall t\in\N,\;\forall x\in\sass^n;
\end{equation}
\item
$e_k(\arcs)(x_1,\dots,x_n)=\{x_k\}$ whenever the states of agent~$k$
and agents~$i\in\neighbors(k,\arcs)$ are all equal and
$e_k(\arcs)(x_1,\dots,x_n)$ is contained in the relative interior of
the convex hull of the states of agent~$k$ and
agents~$i\in\neighbors(k,\arcs)$ whenever the states of agent~$k$ and
agents~$i\in\neighbors(k,\arcs)$ are not all equal.
\end{enumerate}
\end{assumption}\noindent
\begin{remark}
Recall that the values of a set-valued function (from Euclidean space
to Euclidean space) are, by definition, closed
sets. Assumption~\ref{a:f2} thus implies that
$e_k(\arcs)(x_1,\dots,x_n)$ is a compact set contained in the relative
interior of the convex hull of the states of agent~$k$ and
agents~$i\in\neighbors(k,\arcs)$ whenever the states of agent~$k$ and
agents~$i\in\neighbors(k,\arcs)$ are not all equal.
\end{remark}

Assumption~\ref{a:f1} captures the constraints that the
communication topology is supposed to impose on the dynamics. 
Assumption~\ref{a:f2} is a strict convexity assumption.  It implies
that the state~$x_k(t+1)$ is a strict convex combination of the
state~$x_k(t)$ and the states~$x_i(t)$ of
agents~$i\in\neighbors(k,\arcs(t))$.
Assumptions~\ref{a:f1} and~\ref{a:f2} are satisfied in various
examples reported in the literature.
\begin{example}[linear]
The linear model of Section~\ref{s:linear} falls within the scope of
the present study.  Consider a sequence of weighted directed
graphs~$\GG(t)=(\nodes,\arcs(t),w(t))$ with common vertex
set~$\nodes=\{1,\dots,n\}$ and where~$t\in\N$.  Assume the existence
of~$0<e_\mathrm{min}\leq e_\mathrm{max}<\infty$ such that the weight
functions~$w(t)$ take values in~$[e_\mathrm{min},\,e_\mathrm{max}]$
for all~$t\in\N$.  The linear model studied in Section~\ref{s:linear}
corresponds to the case of~$X=\R$ and $f(t,x)=A(\GG(t))x$ with the
matrix~$A(\GG)$~defined in Eq.~\eqref{e:a}.  The discrete-time map~$f$
satisfies Assumptions~\ref{a:f1}--\ref{a:f2} with $e(\arcs)(x)$ being
the set of all possible values~$A(\GG)x$ that are obtained by
considering all possible weighted directed graphs~$\GG$ with arc
set~$\arcs$ and with weights contained in the compact
interval~$[e_\mathrm{min},\,e_\mathrm{max}]$.
\end{example}
\begin{example}[synchronization]\label{e:synchronization}
Consider a population of oscillators that share a common state
space~$S^1$.  Denote the state of oscillator~$k$ by~$\theta_k$ and
consider for each communication graph~$(\nodes,\arcs)$ the Kuramoto
equation~\cite{Ku:84,St:00}
\begin{equation}
\dot{\theta}_k=\omega_k+\sum_{i\in\neighbors(k,\arcs)}\sin(\theta_i-\theta_k),\qquad k=1,\dots,n,
\end{equation}
where we have taken the coupling strenth equal to~$1$.
Let us consider the case that the~$\omega_k$ are all equal.  In this
case we may assume without loss of generality (by introducing a
rotating reference frame) that~$\omega_1=\dots=\omega_n=0$.  We
restrict attention to angles in the interval~$(-\pi/2,\,\pi/2)$ and
introduce local
coordinates~$(-\pi/2,\,\pi/2)\rightarrow\R:\theta_k\mapsto
x_k=\tan(\theta_k)$. Expressed in terms of these local coordinates,
the Kuramoto equation becomes
\begin{equation}\label{e:Kura2}
\dot{x}_k=\sum_{i\in\neighbors(k,\arcs)}\frac{x_i-x_k}{\sqrt{1+x_i^2}\sqrt{1+x_k^2}},\qquad k=1,\dots,n.
\end{equation}
The continuous-time system~\eqref{e:Kura2} may be studied within the
present framework by introducing its
time-$1$~map~$\psi(\arcs):\R^n\rightarrow\R^n$.
Explicitly, consider a time-sequence~$(\nodes,\arcs(t))$ ($t\in\N$)
and put $f(t,x)=\psi(\arcs(t))(x)$.  This discrete-time map~$f$ satisfies
Assumptions~\ref{a:f1} and~\ref{a:f2}
with~$e(\arcs)(x)=\{\psi(\arcs)(x)\}$.  Whereas most papers on
synchronization restrict attention to time-independent, bidirectional
coupling, the present paper allows to study synchronization of
oscillators with time-varying and unidirectional coupling.
\end{example}\noindent
\begin{example}[consensus algorithm]
The following nonlinear consensus protocol is studied in~\cite{OlMu:acc03}.
Let~$\sass=\R$ and consider for a given communication
graph~$(\nodes,\arcs)$ the continuous-time system on~$X^n$
\begin{equation}\label{e:OlMu}
\dot{x}_k=\sum_{i\in\neighbors(k,\arcs)}\gamma_{ik}(x_i-x_k),\qquad k=1,\dots,n.
\end{equation}
where the functions~$\gamma_{ik}:\R\rightarrow\R$ are uneven, locally
Lipschitz and strictly increasing for all $(i,k)\in\arcs$.  
The continuous-time system~\eqref{e:Kura2} may be studied within the
present framework by introducing its
time-$1$~map~$\psi(\arcs):\R^n\rightarrow\R^n$.
Explicitly, consider a time-sequence~$(\nodes,\arcs(t))$ ($t\in\N$)
and put $f(t,x)=\psi(\arcs(t))(x)$.  This discrete-time map~$f$
satisfies Assumptions~\ref{a:f1} and~\ref{a:f2}
with~$e(\arcs)(x)=\{\psi(\arcs)(x)\}$.  Whereas \cite{OlMu:acc03}
restricts attention to time-independent and bidirectional
communication graphs, the present paper allows to study nonlinear
consensus protocols with time-varying and unidirectional
communication.
\end{example}\noindent
\begin{example}[swarming]
The paper~\cite{ViCzJaCoSc:95} proposes a simple model to investigate
self-ordered motion in systems of particles with biologically
motivated interaction.  Their model consists of a collection of
particles moving with constant velocity~$v$. At each time step a given
particle assumes the average direction of motion of the particles in
its neighborhood with some random perturbation added. Denote the
direction of motion of particle~$k$ by~$\theta_k\in S^1$.
The following model is studied in~\cite{ViCzJaCoSc:95}
\begin{equation}\label{e:Vi}
\theta_k(t+1)=\arctan\left(
\frac{\sum_{i\in\{k\}\cup\neighbors(k,\arcs(t))}\sin(\theta_i(t))}%
{\sum_{i\in\{k\}\cup\neighbors(k,\arcs(t))}\cos(\theta_i(t))}
\right),\qquad k=1,\dots,n,
\end{equation}
where we have omitted an additive disturbance term.  Here~$\arcs(t)$
characterizes the nearest neighbor coupling at time~$t$.  We
restrict attention to angles in the interval~$(-\pi/2,\,\pi/2)$ and
introduce local
coordinates~$(-\pi/2,\,\pi/2)\rightarrow\R:\theta_k\mapsto
x_k=\tan(\theta_k)$, as in Example~\ref{e:synchronization}. The
model~\eqref{e:Vi} expressed in terms of these local coordinates falls
within the scope of the present paper. Whereas the study
in~\cite{ViCzJaCoSc:95} is based solely on simulations, the present
appraoch enables us to make analytical statements about convergence.

If the headings~$\theta_k$ are close to each other, then Eq.~\eqref{e:Vi}
may be approximated by the linear equation
\begin{equation}\label{e:Ja}
x_k(t+1)=\sum_{i\in\{k\}\cup\neighbors(k,\arcs(t))}x_i(t),\qquad k=1,\dots,n.%
\end{equation}
with~$x_k\in\R$. This is the swarming model studied
in~\cite{JaLiMo:03} and corresponds to the linear example of
Section~\ref{s:linear} with all weights equal to~$1$.  Hence
Eq.~\ref{e:Ja} also falls within the scope of the present paper.
\end{example}\noindent

Assumption~\ref{a:f2} plays a central role in the forthcoming
analysis. Its importance is already reflected in the following simple
but appealing result.
\begin{lemma}\label{l:convex}
Consider the Data~\ref{d:a}--\ref{d:f} satisfying
Assumption~\ref{a:f2}.  The set-valued
function~$\widetilde{V}:\sass^n\setarrow \sass$ defined for each
$x=(x_1,\dots,x_n)\in \sass^n$ as the convex hull
of~$\{x_1,\dots,x_n\}$, denoted by~$\convexhull\{x_1,\dots,x_n\}$,
is non-increasing along the solutions of Eq.~\eqref{e:consensus};
that is,
\begin{equation}\label{e:nonincrease}
\widetilde{V}(f(t,x))\subseteq \widetilde{V}(x)\qquad \forall t\in\N,\,\forall x\in \sass^n.
\end{equation}
\end{lemma}\noindent
\begin{proof}
Assumption~\ref{a:f2} implies
\begin{equation}
f_k(t,x)\in e_k(\arcs(t))(x)\subseteq\convexhull\{x_1,\dots,x_n\},\qquad
\forall k\in\{1,\dots,n\},\;\forall t\in\N,\;\forall x\in\sass^n,
\end{equation}
from which Eq.~\eqref{e:nonincrease} follows immediately.
\end{proof}\noindent
Within the current framework, the statement of Lemma~\ref{l:convex}
may seem almost trivial; it arises as a straightforward consequence of
Assumption~\ref{a:f2}.  However, the nontrivial contribution lies
precisely in the observation that many examples reported in the
literature satisfy Assumption~\ref{a:f2} and hence may be studied by
means of Lemma~\ref{l:convex}.  The set-valued
function~$\widetilde{V}$ serves as a measure for disagreement.  In the
forthcoming sections, we will use (a slight modification
of)~$\widetilde{V}$ as a set-valued Lyapunov function. The set-valued
nature of the Lyapunov function turns out to be convenient as it
enables the application of Lyapunov techniques in the present context,
where convergence to non-isolated equilibria will be shown.

We end the section with some remarks about the limitations of the
present model.
\begin{remark}[non-strict convexity]
Assumption~\ref{a:f2} is a strict convexity assumption.  For
some applications, it may be desirable to relax this assumption. For
example, the discrete-time system
\begin{equation}
{x}_k(t+1)=\max_{i\in\{k\}\cup\neighbors(k,\arcs(t))}x_i(t),\qquad k=1,\dots,n
\end{equation}
does not satisfy Assumption~\ref{a:f2}, but satisfies instead a
non-strict convexity assumption: $x_k(t+1)$ belongs to (the boundary
of) the convex hull of the state~$x_k(t)$ and the states~$x_i(t)$ of
agents~$i\in\neighbors(k,\arcs(t))$.  The study of dynamical equations
satisfying a non-strict convexity assumption instead of
Assumption~\ref{a:f2} is beyond the scope of the present paper.
\end{remark}\noindent
\begin{remark}[non-Euclidean space]
By assuming that the common state space for the individual agents is
Euclidean, we exclude some global phenomena arising, for example, in
synchronization of periodic motions and in attitude
allignment problems~\cite{SmHaLe:01}.  Indeed, the natural state space
for the individual agents in these examples is respectively $S^1$ and
$SO(3)$ and the global features of the dynamics that are related to
the nontrivial topology of these manifolds, fall outside the scope of
the present method.  Nevertheless, local issues can often be studied
within the present framework, for example, by introducing suitable
coordinate charts. This has already been illustrated in
Example~\ref{e:synchronization}.
\end{remark}\noindent

\section{Stability definitions}
\label{s:stab}

In order to enable a clear and precise formulation of the stability
and convergence properties of the discrete-time
system~\eqref{e:consensus}, we extend the familiar stability concepts
of Lyapunov theory to the present framework.  Notice that we are
interested in the agents' states converging to a common, constant
value and that we expect this common value to depend continuously on
the agents' initial states.  This means that the classical stability
concepts developed for the study of individual (typically isolated)
equilibria are not well-adapted to the present situation.
Alternatively, one may shift attention away from the individual
equilibria and consider the stability properties of the set of
equilibria.  However, set stability does not fully capture the
convergence properties that we are aiming at.  This is illustrated,
for example, by the well-known phenomenon that a trajectory may
converge to the set of equilibria without converging to any of the
individual equilibria.  The stability notions that we introduce below
incorporate, among others, the requirement that all
trajectories converge to one of the equilibria.  We believe that the
notions introduced here constitute a natural framework for questions
studied in coordinated control, formation stabilization,
synchronization, etc.

In the following definition we make a conceptual distinction between
equilibrium solutions and equilibrium points: an equilibrium point is
an element of the state space which is the constant value of an
equilibrium solution.  By referring explicitly to equilibrium
solutions in the following definition, we distinguish the present
stability concepts from the more familiar set stability concepts.
\begin{definition}[stability]
\label{d:stab}
Let~$\mass$ be a finite-dimensional Euclidean space and consider a
continuous map~$f:\N\times \mass\rightarrow\mass$ giving rise to the
discrete-time system
\begin{equation}
x(t+1)=f(t,x(t)).\label{e:generaldynamics}
\end{equation}
Consider a collection of equilibrium solutions of
Eq.~\eqref{e:generaldynamics} and
denote the corresponding set of equilibrium points by~$\Phi$.  With
respect to the considered collection of equilibrium solutions, the
dynamical system~\eqref{e:generaldynamics} is called
\begin{enumerate}
\item\label{d:s}
{\em stable} if 
for each $\phi_1\in\Phi$, 
for all $c_2>0$ and 
for all $t_0\in\N$ 
there is $c_1>0$ such that 
every solution~$\zeta$ of Eq.~\eqref{e:generaldynamics} satisfies: 
if $|\zeta(t_0)-\phi_1|<c_1$
then there is $\phi_2\in\Phi$ such that 
$|\zeta(t)-\phi_2|<c_2$ for all $t\geq t_0$;
\item\label{d:b}
{\em bounded} if
for each $\phi_1\in\Phi$, 
for all $c_1>0$ and 
for all $t_0\in\N$ 
there is $c_2>0$ such that 
every solution~$\zeta$ of Eq.~\eqref{e:generaldynamics} satisfies: 
if $|\zeta(t_0)-\phi_1|<c_1$
then there is $\phi_2\in\Phi$ such that 
$|\zeta(t)-\phi_2|<c_2$ for all $t\geq t_0$;
\item\label{d:ga}
{\em globally attractive} if 
for each $\phi_1\in\Phi$, 
for all $c_1,c_2>0$ and 
for all $t_0\in\N$ 
there is $T\geq 0$ such that 
every solution~$\zeta$ of Eq.~\eqref{e:generaldynamics} satisfies: 
if $|\zeta(t_0)-\phi_1|<c_1$ 
then there is $\phi_2\in\Phi$ such that
$|\zeta(t)-\phi_2|<c_2$ for all $t\geq t_0+T$;
\item\label{d:gas}
{\em globally asymptotically stable} if
it is stable, bounded and globally attractive.
\end{enumerate}
\end{definition}\noindent
\begin{remark}
We may take without loss of generality~$\phi_2=\phi_1$ in
Items~\ref{d:s} and~\ref{d:b} of Definition~\ref{d:stab} (but not in
Item~\ref{d:ga}). Compare with the proof of Theorem~\ref{t:lyapunov}
below.  We have chosen not to incorporate~$\phi_2=\phi_1$ in
Items~\ref{d:s} and~\ref{d:b} in order to make the similarities with
Item~\ref{d:ga} more transparent.
\end{remark}\noindent
Definition~\ref{d:stab} may be interpreted as follows. Stability and
boundedness require that any solution of Eq.~\eqref{e:generaldynamics}
which is initially close to~$\Phi$ remains close to {\em one of the
equilibria in~$\Phi$}, thus excluding, for example, the possibility of
drift along the set~$\Phi$.  Global attractivity implies that every
solution of Eq.~\eqref{e:generaldynamics} converges to one of the
equilibria in~$\Phi$.  If the collection of equilibrium solutions is a
singleton consisting of one equilibrium solution, then the notions of
stability, boundedness, global attractivity and global asymptotic
stability of Definition~\ref{d:stab} coincide with the classical
notions that have been introduced for the study of individual
equilibria.

Slightly stronger stability notions result when uniformity with
respect to initial time is introduced in Definition~\ref{d:stab}.
If the number~$c_1$ (respectively~$c_2$ and~$T$) may be chosen
independently of~$t_0$ in Item~\ref{d:s} (respectively Items~\ref{d:b}
and~\ref{d:ga}) then the dynamical
system~\eqref{e:generaldynamics} is called uniformly stable
(respectively uniformly bounded and uniformly globally attractive)
with respect to the considered collection of equilibrium solutions.

Theorem~\ref{t:lyapunov} below provides a sufficient condition for
uniform stability, uniform boundedness and uniform global asymptotic
stability in terms of the existence of a set-valued Lyapunov function.
This result is convenient since, on the one hand, set-valued Lyapunov
functions arise naturally within the present context as illustrated by
Lemma~\ref{l:convex}, and on the other hand, the stability notions
that are asserted in Theorem~\ref{t:lyapunov} are precisely those that
we are aiming to prove. Apart from its application in the present
paper, Theorem~\ref{t:lyapunov} may be of independent interest.
\begin{theorem}[Lyapunov characterization]\label{t:lyapunov}
Let~$\mass$ be a finite-dimensional Euclidean space and consider a
continuous map~$f:\N\times \mass\rightarrow\mass$ giving rise to the
discrete-time system~\eqref{e:generaldynamics}.  Let~$\Xi$ be a
collection of equilibrium solutions of Eq.~\eqref{e:generaldynamics}
and denote the corresponding set of equilibrium points by~$\Phi$.
Consider
an upper semicontinuous set-valued function ${V}:\mass\setarrow
\mass$ satisfying
\begin{enumerate}
\item\label{a:1}
\(
x\in V(x),\qquad\forall x\in \mass;
\)
\item\label{a:2}
\(
V(f(t,x))\subseteq V(x),\qquad\forall t\in\N,\;\forall x\in \mass.
\)
\end{enumerate}
If ${V}(\phi)=\{\phi\}$ for all~$\phi\in\Phi$ then the dynamical
system~\eqref{e:generaldynamics} is uniformly stable with respect
to~$\Xi$.  If ${V}(x)$ is bounded for all~$x\in \mass$ then the
dynamical system~\eqref{e:generaldynamics} is uniformly bounded with
respect to~$\Xi$.

Consider in addition a
function~$\mu:\image(V)\rightarrow\R_{\geq 0}$ and a lower
semicontinuous function~$\beta:\mass\rightarrow\R_{\geq 0}$ satisfying
\begin{enumerate}
\setcounter{enumi}{2}
\item\label{a:3}
$\mu\circ V: \mass\rightarrow\R_{\geq
0}:x\mapsto\mu(V(x))$ is bounded uniformly with respect
to~$x$ in bounded subsets of~$\mass$;
\item\label{a:4} $\beta$ is positive definite with respect to~$\Phi$;
that is, $\beta(\phi)=0$ for all $\phi\in\Phi$ and $\beta(x)>0$ for
all $x\in \mass\setminus\Phi$;
\item\label{a:5}
$\mu(V(f(t,x)))-\mu(V(x))\leq-\beta(x),\qquad\forall t\in\N,\;\forall x\in \mass$.
\end{enumerate}
If ${V}(\phi)=\{\phi\}$ for all~$\phi\in\Phi$ and ${V}(x)$ is
bounded for all~$x\in \mass$ then the dynamical
system~\eqref{e:generaldynamics} is uniformly globally asymptotically
stable with respect to~$\Xi$.
\end{theorem}\noindent
We briefly comment on the role of the functions~$V$, $\mu$
and~$\beta$ in Theorem~\ref{t:lyapunov}.
The set-valued function~$V$ plays the role of a Lyapunov function
which is non-increasing (decreasing) along the solutions of
Eq.~\eqref{e:generaldynamics}.  The set-valued nature of~$V$ is
crucial. A set-valued function allows for a continuum of minimal
elements which are {\em{not}} comparable with each other.  For this
reason, a set-valued Lyapunov function, unlike a real Lyapunov
function, may be used to conclude that each trajectory converges to
one equilibrium out of a continuum of equilibria.
The function~$\mu$ serves as a measure for the size of the values
of~$V$.  In the present paper, we let~$\mu(V(x))$ be the diameter of
the set~$V(x)$.  The function~$\beta$ characterizes the decrease
of~$V$ along the solutions of Eq.~\eqref{e:generaldynamics} as
measured in terms of~$\mu$.

\begin{proof}[Proof of Theorem~\ref{t:lyapunov}.]
(Uniform stability.) 
Consider arbitrary~$\phi_1\in\Phi$ and~$c_2>0$.
If~${V}(\phi_1)=\{\phi_1\}$ then, by upper semicontinuity of~${V}$, there is~$c_1>0$ such
that~${V}(x)\subset\ball(\phi_1,c_2)$ for
all~$x\in\ball(\phi_1,c_1)$.
Consider arbitrary~$t_0\in\N$ and~$x_0\in\ball(\phi_1,c_1)$ and
let~$\zeta$ denote the solution of Eq.~\eqref{e:generaldynamics} with
$\zeta(t_0)=x_0$. Conditions~\ref{a:1} and~\ref{a:2} of the theorem imply that
\begin{equation}
\zeta(t)\in V(\zeta(t))\subseteq V(x_0)\subset\ball(\phi_1,c_2),\qquad\forall t\geq t_0.
\end{equation}

(Uniform boundedness.)
Consider arbitrary~$\phi_1\in\Phi$ and~$c_1>0$.  If~${V}(x)$ is
bounded for all~$x\in \mass$ then, by upper semicontinuity of~${V}$,
there is~$c_2>0$ such that~${V}(x)\subset\ball(\phi_1,c_2)$ for
all~$x\in\ball(\phi_1,c_1)$.  Consider arbitrary~$t_0\in\N$
and~$x_0\in\ball(\phi_1,c_1)$ and let~$\zeta$ denote the solution of
Eq.~\eqref{e:generaldynamics} with $\zeta(t_0)=x_0$.
Conditions~\ref{a:1} and~\ref{a:2} of the theorem imply that
\begin{equation}
\zeta(t)\in V(\zeta(t))\subseteq V(x_0)\subset\ball(\phi_1,c_2),\qquad\forall t\geq t_0.
\end{equation}

(Uniform global asymptotic stability.)
We have already proved that, if~${V}(\phi)=\{\phi\}$ for
all~$\phi\in\Phi$ and ${V}(x)$ is bounded for all~$x\in \mass$,
then the dynamical system~\eqref{e:generaldynamics} is uniformly
stable and uniformly bounded with respect to~$\Xi$. It remains to
prove uniform global attractivity with respect to~$\Xi$.

Consider arbitrary~$\phi_1\in\Phi$ and~$c_1>0$.
If~${V}(x)$ is
bounded for all~$x\in \mass$ then, by upper semicontinuity of~${V}$,
there is a compact set~$K\subset\mass$ such that~${V}(x)\subseteq K$
for all~$x\in\ball(\phi_1,c_1)$.  
Similarly as above, Conditions~\ref{a:1} and~\ref{a:2} of the theorem imply that
every solution of Eq.~\eqref{e:generaldynamics} initiated in~$\ball(\phi_1,c_1)$ remains in $K$.

Consider in addition arbitrary~$c_2>0$.
If~${V}(\phi)=\{\phi\}$ for all~$\phi\in\Phi$ then, by upper
semicontinuity of~${V}$, there is~$c_3>0$ such that for
all~$x\in\ball(\Phi\cap K,c_3)$ there is~$\phi_2\in\Phi$ such
that ${V}(x)\subset\ball(\phi_2,c_2)$.  
Similarly as above,
Conditions~\ref{a:1} and~\ref{a:2} of the theorem imply that every solution of
Eq.~\eqref{e:generaldynamics} entering~$\ball(\Phi\cap K,c_3)$ remains
in a $c_2$-ball around some equilibrium point~$\phi_2\in\Phi$.

It remains to prove the existence of~$T\geq 0$ such that every
solution of Eq.~\eqref{e:generaldynamics} starting in~$\ball(\phi_1,c_1)$ cannot remain longer than
$T$~subsequent times in~$K$ without entering~$\ball(\Phi\cap K,c_3)$.
In agreement with Conditions~\ref{a:3} and~\ref{a:4} of the theorem
and the lower semicontinuity of~$\beta$, we introduce two real numbers
\[
M=\sup_{x\in{\ball(\phi_1,c_1)}}\mu(V(x))<\infty
\]
and
\[
\Delta=\min_{x\in K\setminus\ball(\Phi,c_3)}\beta(x)>0.
\]
Let~$T\geq 0$ be such that $T\Delta>M$.  Consider arbitrary~$t_0\in\N$
and~$x_0\in\ball(\phi_1,c_1)$ and let~$\zeta$ denote the solution of
Eq.~\eqref{e:generaldynamics} with $\zeta(t_0)=x_0$.  Then
Condition~\ref{a:5} of the theorem implies that for some~$t_1\in [t_0,\,t_0+T]$
\[
\zeta(t_1)\in\ball(\Phi\cap K,c_3),
\]
since otherwise $\mu(V(\zeta(t_0+T)))$ would be smaller than
zero, contradicting that~$\mu$ takes only non-negative values.
Putting everything together, we
conclude that for some~$\phi_2\in\Phi$
\begin{equation}
\zeta(t)
\in       V(\zeta(t))
\subseteq V(\zeta(t_1))
\subset   \ball(\phi_2,c_2),
\qquad\forall t\geq t_0+T.
\end{equation}
\end{proof}\noindent
\begin{remark}\label{r:stepdecrease}
The strict decrease condition (Condition~\ref{a:5}) of
Theorem~\ref{t:lyapunov} may be considerably relaxed.  Consider, for
example,
the following condition which requires that~$\mu\circ V$
decreases over time-intervals of length~$\tau$
\begin{enumerate}
\setcounter{enumi}{5}
\item\label{a:5bis}
There is a time~$\tau\in\N$
such that
\begin{equation}
\mu(V(f(t+\tau-1,\dots f(t+1,f(t,x)) \dots)))
-\mu(V(x))\leq-\beta(x),
\qquad\forall t\in\N,\;\forall x\in \mass.
\end{equation}
\end{enumerate}
Theorem~\ref{t:lyapunov} is still true if Condition~\ref{a:5} is
replaced by Condition~\ref{a:5bis}.
\end{remark}
\begin{remark}\label{r:abstract}
The Lyapunov functions featuring in Theorem~\ref{t:lyapunov} are
set-valued functions~$V:\mass\setarrow\mass$, or equivalently,
single-valued
functions~$V:\mass\rightarrow\closedsubsets(\mass)$.  One may
be interested in generalizations of Theorem~\ref{t:lyapunov}
considering more abstract Lyapunov
functions~$V:\mass\rightarrow Z$ with~$Z$ not necessarily
equal to~$\closedsubsets(\mass)$. The determination of appropriate
structures and properties for the set~$Z$ and the function~$V$
enabling Lyapunov-type of results is beyond the scope of this
study. Nevertheless, we point out that a partial ordering of the
set~$Z$ seems to be an essential ingredient in order to be able to
%
%
conclude that every solution converges to one equilibrium out of a
continuum of equilibria. It turns out that the class of set-valued
functions~$V:\mass\setarrow\mass$ is universal in this
context, in the following sense.
Consider a set~$Z$ and let~$\prec$ by a strict partial ordering of~$Z$.%
\footnote{We adhere to the convention of~\cite{Ro:68} according to
which a partial ordering~$\prec$ of a set~$Z$ is a transitive and
antisymmetric relation on~$Z$.  If we never have $z\prec z$ then
$\prec$ is called a strict partial ordering.}
Consider a set~$\mass$ and a
function~$\widetilde{V}:\mass\rightarrow Z$.  Introduce for
every~$x\in\mass$ the subset~$V(x)\subseteq\mass$
defined by
\[
{V}(x)=\{\hat{x}\in\mass:\widetilde{V}(\hat{x})\prec \widetilde{V}(x)\mbox{
or } \widetilde{V}(\hat{x})=\widetilde{V}(x)\}.
\]
Then we have
${V}(x)={V}(\overline{x})$ whenever
$\widetilde{V}(x)=\widetilde{V}(\overline{x})$ and
${V}(x)\subset{V}(\overline{x})$ whenever $\widetilde{V}(x)\prec
\widetilde{V}(\overline{x})$.
The proof of this statement is elementary and therefore omitted.  Not
only does this result provide some motivation for restricting
attention to set-valued maps~$V:\mass\setarrow\mass$ in
Theorem~\ref{t:lyapunov}, it also points towards a constructive method
for finding such set-valued Lyapunov functions.  See, for example,
the proof of Theorem~\ref{t:usub}.
\end{remark}\noindent

We have the following consequence of
Theorem~\ref{t:lyapunov} and Lemma~\ref{l:convex}.
\begin{theorem}\label{t:usub}
Consider the Data~\ref{d:a}--\ref{d:f} satisfying
Assumption~\ref{a:f2}.  The
discrete-time system~\eqref{e:consensus} is uniformly stable and
uniformly bounded with respect to the collection of equilibrium
solutions~$x_1(t)\equiv\dots\equiv x_n(t)\equiv\mathrm{constant}$.
\end{theorem}\noindent
\begin{proof}[Proof of Theorem~\ref{t:usub}]
In order to apply Theorem~\ref{t:lyapunov} we introduce the set-valued
function
${V}:{\sass^n}\setarrow {\sass^n}$ according to
\[
V(x_1,\dots,x_n)
=(\convexhull\{x_1,\dots,x_n\})^n,\qquad\forall x\in{\sass^n}.
\]
The set-valued function~$V$ is derived from~$\widetilde{V}$ following
the procedure of Remark~\ref{r:abstract}
with~$\closedsubsets(\sass),\subset$ playing the role of~$Z,\prec$.
The set-valued function~$V$ is easily seen to be
globally Lipschitz.  Theorem~\ref{t:usub} follows immediately
from Theorem~\ref{t:lyapunov} and Lemma~\ref{l:convex}.
\end{proof}\noindent

\section{Uniform global attractivity}
\label{s:uga}

Contrary to what might be expected from Lemma~\ref{l:convex} or
Theorem~\ref{t:usub}, the dynamics of Eq.~\eqref{e:consensus} may be
surprisingly complex. For example, the agents' states may fail to
converge to a common value, even in the presence of inter-agent
communication.  In the remainder of the paper, a blend of
graph-theoretic and system-theoretic tools will be used in order to
establish necessary and/or sufficient conditions for global
attractivity.

In the present section we restrict attention to global attractivity
uniform with respect to initial time.  We provide a necessary and
sufficient condition for uniform global attractivity of
Eq.~\eqref{e:consensus}. The condition that we present does not
involve the actual discrete-time map~$f$; it only involves the
sequence of communication graphs~$(\nodes,\arcs(t))$.
\begin{theorem}[uniform global attractivity]\label{t:NS-UGA}
Consider the Data~\ref{d:a}--\ref{d:f} satisfying
Assumptions~\ref{a:f1}--\ref{a:f2}.  The discrete-time
system~\eqref{e:consensus} is uniformly globally attractive with
respect to the collection of equilibrium solutions
$x_1(t)\equiv\dots\equiv x_n(t)\equiv\mathrm{constant}$ if and only if
there is $T\geq 0$ such that for all $t_0\in\N$ the sequence of
communication graphs is weakly connected across~$[t_0,\,t_0+T]$.
\end{theorem}\noindent
\begin{proof}[Proof of Theorem~\ref{t:NS-UGA}]
(Only if, proof by contraposition.)  
Assume that for every $T\geq 0$ there is $t_0\in\N$ such that the
sequence of communication graphs is not weakly connected
across~$[t_0,\,t_0+T]$. By Theorem~\ref{t:equivalence} in the
Appendix, this implies that for every $T\geq 0$ there is $t_0\in\N$
and there are nonempty, disjoint subsets~$\LL_1,\LL_2\subset\nodes$
such that $\neighbors(\LL_1,\arcs(t))$ and
$\neighbors(\LL_2,\arcs(t))$ are both empty for all
$t\in[t_0,\,t_0+T]$.  Pick two different elements~$y,\overline{y}\in
\sass$ and consider any solution~$\zeta$ of Eq.~\eqref{e:consensus}
with initial data
\begin{equation}
\zeta_j(t_0)
\begin{cases}
=y&\quad\forall j\in\LL_1,\\
=\overline{y}&\quad\forall j\in\LL_2,\\
\in\convexhull\{y,\overline{y}\}&\quad\forall j\in\nodes\setminus(\LL_1\cup\LL_2).
\end{cases}
\end{equation}
Assumption~\ref{a:f2} implies that, at time~$t_0+T+1$, we
still have
\begin{equation}
\zeta_j(t_0+T+1)
\begin{cases}
=y&\quad\forall j\in\LL_1,\\
=\overline{y}&\quad\forall j\in\LL_2,\\
\in\convexhull\{y,\overline{y}\}&\quad\forall j\in\nodes\setminus(\LL_1\cup\LL_2),
\end{cases}
\end{equation}
since $\neighbors(\LL_1,\arcs(t))$ and $\neighbors(\LL_2,\arcs(t))$
are both empty for all $t\in[t_0,\,t_0+T]$.  As the time~$T$ may be
chosen arbitrarily large, it is not difficult to see that this
contradicts uniform global attractivity of Eq.~\eqref{e:consensus}
with respect to the equilibrium solutions $x_1(t)\equiv\dots\equiv
x_n(t)\equiv\mathrm{constant}$.

(If.)  
Let $T\geq 0$ be such that for all $t_0\in\N$ the sequence of
 communication graphs is weakly connected across~$[t_0,\,t_0+T]$.
Consider an arbitrary solution~$\zeta$ of Eq.~\eqref{e:consensus} and
an arbitrary time~$t_0\in\N$ and assume that the $\zeta_j(t_0)$'s
($j\in\nodes$) are not all equal.  We first show that
${\widetilde{V}}(\zeta(t_0+(n-1)(T+1)))$ is strictly contained
in~${\widetilde{V}}(\zeta(t_0))$, where ${\widetilde{V}}$ is the
set-valued function introduced in Lemma~\ref{l:convex}.

In order to show this, we introduce a number of auxiliary functions.
Denote the $m$~vertices ($2\leq m\leq n$) of the
polytope~${\widetilde{V}}(\zeta(t_0))$ by $y_1,\dots,y_m\in \sass$.
Associate to each vertex~$y_i$ the set-valued function
$a_i:[t_0,\,\infty)\setarrow\nodes$ identifying the agents located at
that vertex:
\begin{equation}
a_i(t)=\{j\in\nodes:\zeta_j(t)=y_i\},\qquad\forall t\in [t_0,\,\infty),\,\forall i=1,\dots,m.
\end{equation}
The strict convexity assumption (Assumption~\ref{a:f2}) implies
that an agent which is not located at a vertex~$y_i$ at some time
$t\geq t_0$ can never reach this vertex in finite time:
\begin{equation}
a_i(t+1)\subseteq a_i(t),\qquad\forall t\in [t_0,\,\infty),\,\forall
i=1,\dots,m.
\end{equation}
We now establish a strict decrease property for the functions~$a_i$.
We show that, under the connectivity assumption of the theorem, at
least $m-1$ of the functions $a_1,\dots,a_m$ are empty-valued at
$t=t_0+(n-1)(T+1)$.  Consider any intermediate time
interval~$[t_0+p(T+1),\,t_0+p(T+1)+T]$ ($p=0,\dots,n-2$). The
assumption of the theorem implies the existence of an agent~$k$ which
is connected to all other agents across~$[t_0+p(T+1),\,t_0+p(T+1)+T]$.
We distinguish two cases.
\begin{enumerate}
\item
If agent~$k$ is located in one of the vertices at time~$t_0+p(T+1)$;
that is, if there is~$i\in\{1,\dots,m\}$ such that $k\in
a_i(t_0+p(T+1))$, then all the other
vertices~$j\in\{1,\dots,m\}\setminus\{i\}$ where agents are located at
time $t_0+p(T+1)$ (that is, $a_j(t_0+p(T+1))\not =\emptyset$) satisfy
the following: $\neighbors(a_j(t_0+p(T+1)),\arcs(t))\not =\emptyset$
for some $t\in[t_0+p(T+1),\,t_0+p(T+1)+T]$ and hence, by the strict
convexity assumption (Assumption~\ref{a:f2}),
$a_j(t_0+(p+1)(T+1))$ is strictly contained in $a_j(t_0+p(T+1))$.
\item
If agent~$k$ is not located in one of the vertices at time~$t_0+p(T+1)$; that is,
if $k\notin a_j(t_0+p(T+1))$ for all~$j\in\{1,\dots,m\}$, then
all vertices~$j\in\{1,\dots,m\}$ where agents are located at
time $t_0+p(T+1)$ (that is, $a_j(t_0+p(T+1))\not =\emptyset$) satisfy 
the following: $\neighbors(a_j(t_0+p(T+1)),\arcs(t))\not =\emptyset$
for some $t\in[t_0+p(T+1),\,t_0+p(T+1)+T]$ and hence, by the strict
convexity assumption (Assumption~\ref{a:f2}),
$a_j(t_0+(p+1)(T+1))$ is strictly contained in $a_j(t_0+p(T+1))$.
\end{enumerate}
Since the above argument holds for every time
interval~$[t_0+p(T+1),\,t_0+p(T+1)+T]$ ($p=0,\dots,n-2$) and since the
maximum number of agents at each vertex at time~$t_0$ is not greater
than~$n-1$, we may thus conclude that at least $m-1$ of the functions
$a_1,\dots,a_m$ are empty-valued at $t=t_0+(n-1)(T+1)$.  This
establishes a strong decrease property for the set-valued
function~${\widetilde{V}}$. The
polytope~${\widetilde{V}}(\zeta(t_0+(n-1)(T+1)))$ is strictly
contained in the polytope~${\widetilde{V}}(\zeta(t_0))$ and
both sets have not more than one vertex in common.

In order to prepare for an application of Theorem~\ref{t:lyapunov},
we introduce the set-valued
function $V:{{\sass^n}}\setarrow {\sass^n}$ 
according to
\[
V(x_1,\dots,x_n)
=(\convexhull\{x_1,\dots,x_n\})^n,\qquad\forall t\in\N,\;\forall x\in{\sass^n},
\]
similarly as in the proof of Theorem~\ref{t:usub}.  In addition we
also introduce a real-valued function~$\beta:{\sass^n}\rightarrow\R$
according to
\begin{equation}\label{p:d:beta}
\beta(x)=\inf_{\zeta_0,\dots,\zeta_\tau}\diam(V(\zeta_0))-\diam(V(\zeta_\tau)),
\qquad\tau=(n-1)(T+1)
\end{equation}
where~$\diam$ denotes the diameter of a set and where the infimum is
taken over all sequences~$\zeta_0,\dots,\zeta_\tau$ in~$\sass^n$
satisfying
\begin{equation}\label{e:zeta}
\begin{split}
\zeta_0		&=x,\\
\zeta_1		&\in e(\arcs(t_0))(\zeta_0),\\
       		&\vdots\\
\zeta_\tau	&\in e(\arcs(t_0+\tau-1))(\zeta_{\tau-1}),
\end{split}
\end{equation}
for some~$t_0\in\N$. The~$\zeta_i$ may be interpreted as the states
that are possibly reachable from~$x$ in $i$~time steps according
to Assumption~\ref{a:f2}. 

We establish some
useful properties of~$\beta$.
First, the collection of all sequences~$\zeta_0,\dots,\zeta_\tau$
in~$\sass^n$ satisfying~\eqref{e:zeta} for some~$t_0\in\N$ 
\begin{itemize}
\item
is nonempty and compact for all~$x\in\sass^n$;
\item 
depends continuously on~$x$.
\end{itemize}
This follows from the observations that the set-valued
functions~$e(\arcs):\sass^n\setarrow\sass^n$ are continuous and take
nonempty, compact values (implied by Assumption~\ref{a:f2}) and that
there is only a finite number of possible sequences of
$\tau$~communication graphs.  Secondly, the expression
$\diam(V(\zeta_0))-\diam(V(\zeta_\tau))$
 being minimized
in~\eqref{p:d:beta}
\begin{itemize}
\item 
depends continuously on~$\zeta_0$ and~$\zeta_\tau$: indeed, this
follows from the global Lipschitz character of~$V$;
\item
is zero whenever the $n$~components of~$x$ are all equal: indeed, in
this case Assumption~\ref{a:f2} implies
that~$\zeta_0=\dots=\zeta_\tau=x$ is the only sequence
satisfying~\eqref{e:zeta} for some~$t_0\in\N$;
\item
is strictly positive whenever the $n$~components of~$x$ are not all
equal: indeed, this follows from the strict decrease property
of~$\widetilde{V}$ that has been established above and the observation
that~$V(x)=(\widetilde{V}(x))^n$.
\end{itemize}
Putting all elements together, we conclude that the
function~$\beta:{\sass^n}\rightarrow\R$ is continuous and positive
definite with respect
to~$\{(x_1,\dots,x_n)\in\sass^n:x_1=\dots=x_n\}$.  

The proof is concluded with an application of Theorem~\ref{t:lyapunov}
(see Remark~\ref{r:stepdecrease}) which establishes uniform global
asymptotic stability of Eq.~\eqref{e:consensus}.
\end{proof}

\section{Non-uniform global attractivity}
\label{s:ga}

The previous section has presented a necessary and sufficient
condition for uniform global attractivity.  Not surprisingly, this
necessary and sufficient condition involves a connectivity
requirement on the sequence of directed graphs, uniform with respect
to initial time.  In the present section, we turn attention to global
attractivity, not necessarily uniform with respect to initial time.
Unlike the previous section, where the situation is quite clear, the
study of non-uniform attractivity turns out to be much more subtle.
We start with a necessary condition for global attractivity, not
 necessarily uniform with respect to time.
\begin{proposition}[global attractivity]\label{p:N-GA}
Consider the Data~\ref{d:a}--\ref{d:f} satisfying
Assumptions~\ref{a:f1}--\ref{a:f2}.  Global
attractivity of the discrete-time system~\eqref{e:consensus} with
respect to the collection of equilibrium solutions
$x_1(t)\equiv\dots\equiv x_n(t)\equiv\mathrm{constant}$ implies that
for all $t_0\in\N$ the sequence of communication graphs is weakly
connected across~$[t_0,\,\infty)$.
\end{proposition}\noindent
The proof of Proposition~\ref{p:N-GA} is very similar to the proof of
the {\it{only~if}}-part of Theorem~\ref{t:NS-UGA} and is
omitted.  The necessary condition featuring in
Proposition~\ref{p:N-GA} may be interpreted as a non-uniform version
of the necessary and sufficient condition featuring in
Theorem~\ref{t:NS-UGA}.
It may therefore seem tempting to conjecture that this necessary
condition is also sufficient for global attractivity, not necessarily
uniform with respect to time.  However, the following counterexample
shows that this is not the case.

\subsection{Counterexample}
\label{s:counter}

The counterexample is concerned with three agents sharing a common
state space~$\sass=\R$.
Among all possible directed graphs
on the vertex set~$\{1,2,3\}$ we consider
\begin{equation}
\begin{split}
\arcs_\mathrm{a}&=\{(1,2)\},\\
\arcs_\mathrm{b}&=\{(1,2),(2,1)\},\\
\arcs_\mathrm{c}&=\{(3,2)\},\\
\arcs_\mathrm{d}&=\{(2,3),(3,2)\}.
\end{split}
\label{e:arcsets}
\end{equation}
We will introduce a sequence of
directed graphs which consists of the concatenation of finite
sequences of the form~$B_s$ ($s\in\{0\}\cup\N$):
\[
B_s=
\underset{2s}{\underbrace{\arcs_\mathrm{a},\dots,\arcs_\mathrm{a}}},
\arcs_\mathrm{b},
\underset{2s+1}{\underbrace{\arcs_\mathrm{c},\dots,\arcs_\mathrm{c}}},
\arcs_\mathrm{d}.
\]

\begin{proposition}\label{p:counter}
Let the number of agents be given by~$n=3$ and consider the sequence
of communication graphs~$(\nodes,\arcs(t))$ with $t\in\N$
corresponding to the concatenation \( B_0,B_1,B_2,B_3,\dots \) Let the
common state space for the individual agents be given by~$\sass=\R$
and consider the discrete-time map~$f$ 
corresponding to the linear example of
Section~\ref{s:linear} with all the weights equal to one.  Let~$\zeta$ be
the solution of Eq.~\eqref{e:consensus} with initial
data~$\zeta_1(1)=0$ and $\zeta_2(1)=\zeta_3(1)=1$.  Then the three
components of~$\zeta(t)$ do not converge to a common value as
$t\rightarrow\infty$.
\end{proposition}
\begin{proof}
In order to show that the three components of~$\zeta(t)$ do not
converge to a common value as $t\rightarrow\infty$, we evaluate the
difference~$\zeta_3-\zeta_1$ at the sequence of time-instants \(
t_1,t_2,t_3,\dots\rightarrow\infty \) determined by \(
t_{p+1}-t_{p}=p+1 \) for all~$p\in\N$ and $t_1=2$
(Fig.~\ref{f:counter}).
\begin{figure}[h]
\hspace*{\fill}
\psfig{file=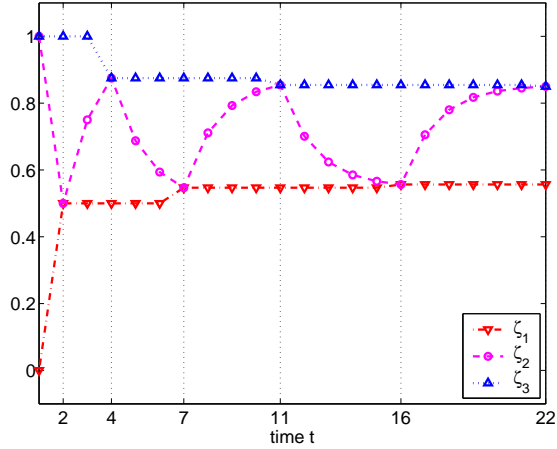,height=6cm}
\hspace*{\fill}
\caption{Counterexample}
\label{f:counter}
\end{figure}
For the ease of notation, let us denote~$\zeta_3(t_p)-\zeta_1(t_p)$
by~$v(p)$.
It is not difficult to see that
\begin{align}
v(p+1)&=v(p)-\frac{1}{2^{p+1}}v(p)=\frac{2^{p+1}-1}{2^{p+1}}v(p),\qquad\forall
p\in\N,\label{e:recursive}\\ v(1)&=1/2,
\end{align}
from which it follows that $0<v(p)<1$ for all~$p\in\N$.
Clearly~$v(p)$ is decreasing with~$p$ and converges
to a limit as~$p\rightarrow\infty$.
The total accumulative decrease of~$v$ satisfies
\begin{align}\label{e:decrease}
\sum_{p=1}^\infty \big( v(p)-v(p+1) \big)
&=\sum_{p=1}^\infty \frac{1}{2^{p+1}}v(p)\\
&<\sum_{p=1}^\infty \frac{1}{2^{p+1}}=\frac{1}{2},
\end{align}
where we have used the recursive relation~\eqref{e:recursive} and the
observation that~$0<v(p)<1$ for all~$p\in\N$.
Since $v(1)=1/2$ we conclude that
\begin{equation}
\lim_{p\rightarrow\infty}\zeta_3(t_p)-\zeta_1(t_p)=
\lim_{p\rightarrow\infty}v(p)=v(1)-\sum_{p=1}^\infty \big( v(p)-v(p+1)
\big) >0.
\end{equation}
\end{proof}

\subsection{Discussion of counterexample}
\begin{enumerate}
\item
First of all, notice that the sequence of directed graphs featuring
in the counterexample satisfies
the necessary condition of Proposition~\ref{p:N-GA}.  The union of the
arc sets~$\arcs(t)$ over any interval of the form~$[t_0,\,\infty)$ is given
by~$\{(1,2),(2,1),(2,3),(3,2)\}$ which corresponds to a connected,
bidirectional graph.  This counterexample thus clearly shows that the
necessary condition of Proposition~\ref{p:N-GA} is not sufficient for
global attractivity.
\item
The counterexample points towards a non-intuitive phenomenon, namely
that more information exchange does not necessarily lead to improved
convergence properties, and may eventually even destroy convergence of
the agents' states to a common value.  To be more precise, if we
replace the arc sets~$\arcs_\mathrm{a}$ and~$\arcs_\mathrm{c}$ by the
empty set, which results in reduced communication, then the
discrete-time system featuring in Proposition~\ref{p:counter} becomes
globally attractive.  By taking~$\arcs_\mathrm{a}$
and~$\arcs_\mathrm{c}$ as defined in~\eqref{e:arcsets} instead of the
empty set, which increases the communication, the agents'
states fail to converge to a common value.  It is, of course,
well-known that more feedback does not necessarily mean improved
convergence, but we are inclined to believe that the present
phenomenon is different from previously observed, similar phenomena.
To illustrate this, let us note that the loss of convergence observed in
the counterexample can only occur in the context of
time-dependent and unidirectional communication (see
Theorem~\ref{t:NS-UGA} and Theorem~\ref{t:bidirectional} of the next
section).
\end{enumerate}

\section{Bidirectional communication}
\label{s:bidirectional}

The study of global attractivity of Eq.~\eqref{e:consensus}, not
necessarily uniform with respect to initial time, is simplified
considerably when bidirectional communication is assumed.
\begin{theorem}[bidirectional]\label{t:bidirectional}
Consider the Data~\ref{d:a}--\ref{d:f} satisfying
Assumptions~\ref{a:f1}--\ref{a:f2}.  Assume in
addition that the communication graphs~$(\nodes,\arcs(t))$ are
bidirectional for all~$t\in\N$.  The discrete-time
system~\eqref{e:consensus} is globally attractive with respect to the
collection of equilibrium solutions $x_1(t)\equiv\dots\equiv
x_n(t)\equiv\mathrm{constant}$ if and only if 
for all~$t_0\in\N$ the sequence of bidirectional
graphs~$(\nodes,\arcs(t))$ is connected across~$[t_0,\,\infty)$.
\end{theorem}\noindent
The proof of Theorem~\ref{t:bidirectional} is based on an analysis
of $\omega$-limit sets. It is, of course, well-known that $\omega$-limit sets
play a central role in the stability analysis of dynamical systems
through the celebrated LaSalle principle, but it is quite remarkable
that the notion of $\omega$-limit set proves to be useful in the
present context, since we are considering time-varying dynamical
systems, not necessarily periodic.
\begin{proof}[Proof of Theorem~\ref{t:bidirectional}]
(Only if.) This follows from Proposition~\ref{p:N-GA} and the
observation that a bidirectional graph is connected if and only if it
is weakly connected.
(If.)
It suffices to prove that every solution of Eq.~\eqref{e:consensus}
converges to one of the equilibrium solutions $x_1(t)\equiv\dots\equiv
x_n(t)\equiv\mathrm{constant}$.  Indeed, by continuity and compactness
arguments, convergence of all individual solutions actually implies
global attractivity in the sense of Definition~\ref{d:stab}.  The
observation that compactness may be invoked is not trivial (since the
set of equilibrium points is unbounded) but follows from the
boundedness property of Eq.~\eqref{e:consensus} established in
Theorem~\ref{t:usub}.

In the remainder of the proof we show that every solution of
Eq.~\eqref{e:consensus} converges to one of the equilibrium solutions
$x_1(t)\equiv\dots\equiv x_n(t)\equiv\mathrm{constant}$.  Consider
arbitrary~$t_0\in\N$ and~$x_0\in{\sass^n}$, and let~$\zeta$ denote the
solution of Eq.~\eqref{e:consensus} with initial
data~$\zeta(t_0)=x_0$.  Let~$\Omega$ denote the $\omega$-limit set
of~$\zeta$.
We start with observing that 
\begin{equation}\label{e:VOmega}
\widetilde{V}(x)=\widetilde{V}(\overline{x}),
\qquad\forall x,\overline{x}\in\Omega,
\end{equation} 
where~$\widetilde{V}$ is the set-valued function introduced in
Lemma~\ref{l:convex}. Indeed, the existence
of~$x,\overline{x}\in\Omega$ with $\widetilde{V}(x)\not
=\widetilde{V}(\overline{x})$ would contradict the non-increase
property established in Lemma~\ref{l:convex}. We denote the constant
value of~$\widetilde{V}$ on~$\Omega$ by~$\widetilde{V}(\Omega)$.

Clearly, in order to establish that~$\zeta$ converges to one of the
equilibrium solutions $x_1(t)\equiv\dots\equiv
x_n(t)\equiv\mathrm{constant}$, it suffices to prove
that~$\widetilde{V}(\Omega)$ is a singleton. We prove this by
contradiction.  Assume that~$\widetilde{V}(\Omega)$ is not a
singleton. In this case, $\widetilde{V}(\Omega)$~is a
polytope with $m$~vertices ($2\leq m\leq n$) which we denote by
$y_1,\dots,y_m\in \sass$.  We associate to each vertex~$y_i$ a
set-valued function $b_i:{\sass^n}\setarrow\nodes$ identifying the agents
located at that vertex:
\begin{equation}
b_i(x)=\{j\in\nodes:x_j=y_i\},\qquad\forall x\in{\sass^n},\,\forall i=1,\dots,m.
\end{equation}
Observation~\eqref{e:VOmega} implies that
\begin{equation}\label{e:VOmega2}
b_i(x)\not =\emptyset,\qquad\forall x\in\Omega,\,\forall i=1,\dots,m.
\end{equation}
We arrive at a contradiction with the aid of the following result.
\begin{claim}\label{claim}
If the connectivity condition mentioned in
Theorem~\ref{t:bidirectional} holds, then for every~$x\in\Omega$ with
$b_i(x)\not =\emptyset$ for all~$i\in\{1,\dots,m\}$ there
exists~$\overline{x}\in\Omega$ with~$b_i(\overline{x})$ strictly
contained in~$b_i(x)$ for some~$i\in\{1,\dots,m\}$.
\end{claim}\noindent
Indeed, a repetitive application of this result eventually leads to
the existence of~$x^*\in\Omega$ with~$b_i(x^*)=\emptyset$ for
some~$i\in\{1,\dots,m\}$, contradicting~\eqref{e:VOmega2}.  We have
thus shown, by contradiction, that~$\widetilde{V}(\Omega)$ is a
singleton and thus that~$\zeta$ converges to one of the equilibrium
solutions $x_1(t)\equiv\dots\equiv x_n(t)\equiv\mathrm{constant}$.

We end the proof of Theorem~\ref{t:bidirectional} with a proof of
Claim~\ref{claim}. Consider an arbitrary~$x\in\Omega$ with $b_i(x)\not
=\emptyset$ for all~$i\in\{1,\dots,m\}$. There is, by definition, a
sequence of times~$t_q$ ($q\in\N$) tending to infinity such
that~$\zeta(t_q)\rightarrow x$ as $q\rightarrow\infty$.
Based upon the sequence of times~$t_q$ we construct another sequence
of times~$\hat{t}_q$ ($q\in\N$) tending to infinity, as follows.  Let
for each~$q\in\N$ the time~$\hat{t}_q$ be defined as the first time
greater than or equal to~$t_q$ at which communication occurs between
an agent~$k\in b_i(x)$ and an agent~$l\not\in b_i(x)$ for
some~$i\in\{1,\dots,m\}$:
\[
\hat{t}_q=\min\{t\in[t_q,\,\infty):\neighbors(b_i(x),\arcs(t))\not
=\emptyset\mbox{ for some~$i=1,\dots,m$}\},\qquad\forall q\in\N.
\]
The validity of this construction follows from the connectivity
condition of Theorem~\ref{t:bidirectional} that we are assuming to
hold.  We may assume without loss of generality (by considering an
appropriate subsequence if necessary) that
\begin{enumerate}
\item the~$\zeta(\hat{t}_q)$ converge to a limit
point~$\hat{x}\in\Omega$ as $q\rightarrow\infty$ (by boundedness of
the solution~$\zeta$);
\item the~$\arcs(\hat{t}_q)$ are all equal, say,
$\arcs(\hat{t}_q)=\widehat{\arcs}$ for all $q\in\N$ (since there is only a
finite number of possible communication graphs);
\item the~$\zeta(\hat{t}_q+1)$ converge to a limit
point~$\overline{x}\in\Omega$ as $q\rightarrow\infty$ (by continuity
of the set-valued map~$e(\widehat{\arcs})$ and compactness
of~$e(\widehat{\arcs})(\hat{x})$ implied by Assumption~\ref{a:f2}).
\end{enumerate}
It follows that~$\overline{x}\in e(\widehat{\arcs})(\hat{x})$.
By construction of the sequence of times~$\hat{t}_q$ we conclude that,
on the one hand, $b_i(\hat{x})=b_i(x)$ for all~$i\in\{1,\dots,m\}$,
and, on the other hand, $b_i(\overline{x})$~is strictly contained
in~$b_i(\hat{x})$ for some~$i\in\{1,\dots,m\}$ (by
Assumption~\ref{a:f2}).  This concludes the proof of claim~\ref{claim}
and hence also of Theorem~\ref{t:bidirectional}.
\end{proof}

\section{Conclusion}
\label{s:conclusion}

We have considered in the paper a simple but appealing model for
$n$~interacting agents via unidirectional and time-dependent
communication. This model finds wide application in a variety of
fields including synchronization, swarming and distributed decision
making.  In the model, each agent updates his state based upon the
information received from other agents, according to a simple rule. By
considering simple dynamics for the individual agents, we are able to
focus on the main aspects of communication topology without having
to deal with the additional complications arising from complex
dynamical agents.

The analysis starts with the assumption that each agent updates his
state according to a strict convex combination of its neighbors'
states, an assumption which is satisfied in various examples studied
in the literature.  This assumption leads to the development of a
set-valued Lyapunov theory, with the convex hull of the individual
agents' states playing the role of a non-increasing Lyapunov function.
Contrary to what might be expected, the dynamics of the multi-agent
system turns out to be quite subtle and a blend of graph-theoretic and
system theoretic tools is used in order to obtain necessary and/or
sufficient conditions for convergence.

The strongest result is obtained for the case of bidirectional
communication.  In that case, it is shown that convergence of the
individual agents' states to a common value is guaranteed if, during each interval of
the form~$[t_0,\infty)$, each agent sends information to each other
agent, either through direct communication or indirectly via
intermediate agents.

The case of unidirectional communication is more subtle.  It is shown
by means of a counterexample that, contrary to what might be expected,
convergence of the individual agents' states to a common value is not
necessarily guaranteed, even if during each interval of the
form~$[t_0,\infty)$ there is an agent who sends information to all
other agents, either through direct communication or indirectly via
intermediate agents.  Convergence is proven, however, if a uniform
bound is imposed on the time it takes for the information to spread
over the network.

The counterexample that is used to show that the agents' states may
fail to converge to a common value, even in the presence of
communication, points towards a counter-intuitive phenomenon.  Namely
that more information exchange does not necessarily lead to improved
convergence and may eventually even lead to a loss of convergence,
even for the simple models studied in the present paper.  The study of
the quantitative relationship between communication topology and speed
of convergence remains an interesting area of further research.


\bibliographystyle{abbrv}
\bibliography{/home/lmoreau/bibliography/literature,/home/lmoreau/bibliography/multi-agent}

\appendix

\section{Graph-theoretic result}
\label{s:appendix}

\begin{theorem}\label{t:equivalence}
A directed graph~$(\nodes,\arcs)$ is weakly connected if and only if
every pair of nonempty, disjoint subsets~$\LL_1,\LL_2\subset\nodes$
satisfies $\neighbors(\LL_1,\arcs)\cup\neighbors(\LL_2,\arcs)\not
=\emptyset$.
\end{theorem}\noindent
\begin{proof}
(Only if.)  Assume that the node~$k$ is connected to all other
nodes, and consider an arbitrary pair of nonempty, disjoint
subsets~$\LL_1,\LL_2\subset\nodes$. 
We distinguish three cases.
\begin{enumerate}
\item If $k\in\LL_1$ then $\neighbors(\LL_2,\arcs)$ is nonempty.
\item If $k\in\LL_2$ then $\neighbors(\LL_1,\arcs)$ is nonempty.
\item If $k\notin\LL_1\cup\LL_2$ then $\neighbors(\LL_1,\arcs)$ and
$\neighbors(\LL_2,\arcs)$ are both nonempty.
\end{enumerate}
In all three cases
$\neighbors(\LL_1,\arcs)\cup\neighbors(\LL_2,\arcs)\not
=\emptyset$.

(If.) The proof consists of a constructive algorithm that is
guaranteed to terminate in a finite number of steps.  Each
step~$i\in\N$ (except the last step) involves the selection of four
nonempty sets~$\LL_1(i)\subseteq\FF_1(i)\subset\nodes$
and~$\LL_2(i)\subseteq\FF_2(i)\subset\nodes$ satisfying:
\begin{itemize}
\item $\FF_1(i)$ and $\FF_2(i)$ are disjoint;
\item each node in $\LL_1(i)$ is connected to each other node in $\FF_1(i)$
and each node in $\LL_2(i)$ is connected to each other node in $\FF_2(i)$.
\end{itemize}

\underline{Step~$1$.\/} Set $\LL_1(1)=\FF_1(1)=\{m_1\}$ and
$\LL_2(1)=\FF_2(1)=\{\overline{m}_1\}$, where $m_1$ and $\overline{m}_1$ are
two arbitrary, different nodes of the graph.  (Here we have assumed
that the graph has at least 2 nodes. If the graph has only one node,
then the statement of Theorem~\ref{t:equivalence} is trivial.)

\underline{Step~$i>1$.\/} As we are proving the {\it{if}}-part of the theorem,
we may assume that
$\neighbors(\LL_1(i-1),\arcs)\cup\neighbors(\LL_2(i-1),\arcs)\not
=\emptyset$. We restrict attention to the case
$\neighbors(\LL_2(i-1),\arcs)\not =\emptyset$ and consider
$m_i\in\neighbors(\LL_2(i-1),\arcs)$.  (The alternative case of
$\neighbors(\LL_1(i-1),\arcs)\not =\emptyset$ can be dealt with in a
completely similar way, be interchanging the role of $\LL_1$ and
$\LL_2$, respectively $\FF_1$ and $\FF_2$.)
We distinguish four cases:
\begin{enumerate}
\item If $m_i\in\FF_1(i-1)$ and $\FF_1(i-1)\cup\FF_2(i-1)=\nodes$, then
the algorithm terminates: any agent $\overline{m}_i\in\LL_1(i-1)$ is
connected to all other agents in the graph.
\item If $m_i\in\FF_1(i-1)$ and $\FF_1(i-1)\cup\FF_2(i-1)\not
=\nodes$, then set
\begin{align}
\LL_1(i)&=\LL_1(i-1),\\
\FF_1(i)&=\FF_1(i-1)\cup\FF_2(i-1),\\
\LL_2(i)=\FF_2(i)&=\{\overline{m}_i\},
\end{align}
where $\overline{m}_i$ is an arbitrary node which does not belong to
$\FF_1(i-1)\cup\FF_2(i-1)$.
\item If $m_i\not\in\FF_1(i-1)\cup\FF_2(i-1)$, then set
\begin{align}
\LL_1(i)&=\LL_1(i-1),\\
\FF_1(i)&=\FF_1(i-1),\\
\LL_2(i)&=\{m_i\},\\
\FF_2(i)&=\FF_2(i-1)\cup\{m_i\}.
\end{align}
\item If $m_i\in\FF_2(i-1)\setminus\LL_2(i-1)$ then
\begin{align}
\LL_1(i)&=\LL_1(i-1),\\
\FF_1(i)&=\FF_1(i-1),\\
\LL_2(i)&=\LL_2(i-1)\cup\{m_i\},\\
\FF_2(i)&=\FF_2(i-1).
\end{align}
\end{enumerate}

The algorithm is guaranteed to end in a finite number of steps, since
at each step (except at the last step) either the number of agents in
$\FF_1\cup\FF_2$ increases, or the number of agents in
$\FF_1\cup\FF_2$ remains unchanged and the number of agents in
$\LL_1\cup\LL_2$ increases.
\end{proof}

\end{document}